\documentclass{amsart}
\usepackage{amsmath,amsfonts,amsthm,enumerate}

\newcommand{\End}{\text{\rm End}}
\newcommand{\thmref}[1]{Theorem~\ref{#1}}

\newcommand{\corref}[1]{Corollary~\ref{#1}}

\newcommand{\eqnref}[1]{~(\ref{#1})}

\newcommand{\germ}{\mathfrak}
\newtheorem{thm}{Theorem}[section]
\newtheorem{lem}[thm]{Lemma}

\theoremstyle{definition}
\newtheorem{cor}[thm]{Corollary}

\newtheorem{example}[thm]{Example}

\newtheorem{rem}[thm]{Remark}

\theoremstyle{rem}
\numberwithin{equation}{section}

\def\D{\Delta}
\def\C{{\mathbb C}}
\def\G{{\mathfrak {g}}}

\def\H{{\mathcal {H}}}
\def\Z{{\mathbb Z}}

\title{}
\begin{document}
\subjclass{Primary 17B67, 81R10}

\title{Structure of intermediate Wakimoto modules}
\author{Ben L.  Cox}
\author{Vyacheslav  Futorny}
\address{Department of Mathematics \\
College of Charleston \\
66 George Street  \\
Charleston SC 29424, USA}\email{coxbl@cofc.edu}

\address{School of Mathematics and Statistics\\
 University of Sydney\\
Sydney 2006, Australia\\
 On leave from Institute of Mathematics\\
University of S\~ao Paulo \\
Caixa Postal 66281  \\
S\~ao Paulo, CEP 05315-970, Brazil}\email{futorny@ime.usp.br}

\begin{abstract}  We show that our construction of boson type
 realizations of
affine $\mathfrak{sl}(n+1)$ in terms of intermediate Wakimoto modules gives representations that are generically isomorphic to certain Verma type modules.   We then use this identification to obtain information about the submodule structure of intermediate Wakimoto modules.
  \end{abstract}


\maketitle
\section{Verma type modules}

Modules induced from the {\it natural Borel subalgebra} were first introduced by H.~Jakobsen and V.~Kac in their study of unitarizable highest weight representations of affine Kac-Moody algebras (see \cite{JK}). They were studied in
\cite{MR95a:17030} under the name of  {\it imaginary Verma modules}.
A Fock space realization of the imaginary Verma modules for
$\hat{\mathfrak{sl}}(2)$ were constructed by D.~Bernard and G.~Felder
in \cite{BF} and then extended in \cite{C2} to the case of
$\hat{\mathfrak{sl}}(n)$. These realizations are given generically by  certain
Wakimoto type modules.

In his effort to prove the Kac-Kazhdan conjecture on the characters of irreducible representations of afffine Kac-Moody algebras at the critical level, M. Wakimoto discovered a remarkable boson realization of $\hat{\mathfrak{sl}}(2)$ on the Fock space $\mathbb C[x_i,y_j\,|\,i\in\mathbb Z, j\in\mathbb N]$ \cite{W}.
Wakimoto modules for general affine Lie algebras were introduced by B.~Feigin and
E.~Frenkel in \cite{MR89k:17016} by a homological characterization, \cite{MR92d:17025} which play an important role in the conformal field theory providing a new bosonization rule for the Wess-Zumino-Witten models.
Wakimoto modules have a geometric interpretation as certain sheaves on a semi-infinite flag manifold \cite{MR92f:17026}. They  belong to the category $\mathcal O$
and generically are isomorphic to corresponding Verma modules.  Explicit formulae for these realizations for $\hat{\mathfrak{sl}}(2)$ are given in  \cite{MR89k:17016} and for general affine Lie algebras they are given in \cite{MR1482938} and  \cite{MR1412381}.

Affine Lie algebras admit Verma type modules associated with
 non-standard Borel subalgebras, see \cite{C}, \cite{FS}
and \cite{MR89m:17032}.
In our work \cite{MR2065632} the problem of finding  suitable
boson type realizations for all Verma type modules over $\hat{\mathfrak{sl}}(n+1)$ was solved.
In Theorem \ref{realization} below we construct
  such realizations, {\it intermediate Wakimoto modules},
 for a series of generic Verma type modules depending on
 the parameter $0\leq r\leq n$. If $r=n$ this construction
coincides with the boson realization of Wakimoto modules in
\cite{MR89k:17016}, and when $r=0$ this is a
 realization described in \cite{C2}. One difficulty that arises in the
study of Verma type modules that are not induced from a standard
Borel subalgebra is that certain of their weight spaces are
infinite dimensional. On the other hand the structure of
representations that have infinite dimensional weight spaces is an
important problem that appears naturally in other contexts.
Besides appearing in the representation theory of infinite
dimensional Heisenberg Lie algebras such representations also
arise in the work of \cite{MR95i:17007}, \cite{MR88h:17022},
\cite{MR2001k:17011}, \cite{MR2003g:17034},  and
\cite{MR2002f:17041}.  Intermediate Wakimoto modules form another
family of representations with certain weight spaces being
infinite dimensional. Using the realization of the intermediate
Wakimoto module for $\hat{\mathfrak{sl}}(n+1, \mathbb C)$ given in
\cite{MR2065632} and below,  we show that generically Verma type
modules and intermediate Wakimoto modules are isomorphic, which is
an analog of the classical relation between Verma and Wakimoto
modules (see \corref{irred} below). Moreover, when intermediate
Wakimoto modules are in general position (but not necessarily
isomorphic to Verma type module) we completely describe their
submodule structure (\ref{cor30}).

Verma type modules have a complicated structure when the center
$c$ acts by zero (see for example \cite{MR95a:17030}). The
realization  given in \thmref{realization} yields information
about the structure of these modules at least in the case of
$\hat{\mathfrak{sl}}(2, \mathbb C)$. In the last section we
present the formulas for the singular elements in Imaginary Verma
modules recently obtained by B.Wilson \cite{Wi}. These formulas
were inspired by the free field realization of Imaginary Verma
modules for $\hat{\mathfrak{sl}}(2,\mathbb C)$.

\section{Verma type modules}

Fix a positive  integer
$n$, $0\leq r\leq n$,
$\gamma\in\mathbb C^*$. Set
$k=\gamma^2-(r+1)$.
Let $\mathfrak g=\mathfrak{sl}(n+1,\mathbb C)$ and let
 $E^{ij}$,  $i,j=1, \ldots, n+1$ be
the standard basis for  $\mathfrak{gl}(n+1,\mathbb C)$.  Set
$H_i:=E^{ii}-E^{i+1,i+1}$, $E_i:=E^{i,i+1}$, $F_i:=E^{i+1,i}$ which is a
basis for $\mathfrak{sl}(n+1,\mathbb C)$. Furthermore
we denote the Killing form by
$(X|Y)=\text{tr}\,(XY)$ and $X_m=t^m\otimes X$ for
$X,Y\in \mathfrak g$ and $m\in\mathbb Z$.  Let
$\{\alpha_1,\dots,\alpha_n\}$ be a base for $\Delta^+$, the positive set
of
roots for $\mathfrak g$, such that $H_i=\check\alpha_i$
and let
$\Delta_r$ be the root system with base $\{\alpha_1,\dots,\alpha_r\}$
($\Delta_r=\emptyset$, if $r=0$) of the Lie subalgebra $\mathfrak g_r=
\mathfrak{sl}(r+1,\mathbb C)$. A Cartan subalgebra
$\mathfrak H$ (respectively $\mathfrak H_r$) of $\mathfrak g$
(respectively $\mathfrak g_r$) is spanned by $H_i$, $i=1, \ldots, n$
(respectively $i=1, \ldots, r$) and set $\mathfrak H_0=0$.

 For any Lie algebra $\mathfrak a$, let
$L(\mathfrak a
)=\mathfrak a \bigotimes\mathbb{C} [t,t^{-1}]$ be the loop algebra of
$\mathfrak a$. Then   $\hat{\mathfrak g}=\hat{\mathfrak{sl}}(n+1,\mathbb C) =L({\mathfrak g}
)\oplus\mathbb{C}
c\oplus\mathbb{C} d$ and $\hat{\mathfrak g}_r=L({\mathfrak g_r}
)\oplus\mathbb{C}
c\oplus\mathbb{C} d$
 are the associated  affine Kac-Moody
algebras with $\hat{\mathfrak H}=\mathfrak H\oplus\mathbb{C}
c\oplus\mathbb{C} d$ and $\hat{\mathfrak H}_r=\mathfrak H_r\oplus\mathbb{C}
c\oplus\mathbb{C} d$ respectively.

The  algebra
$\hat{\mathfrak g}$ has generators
$E_{im}, F_{im}, H_{im}$,  $i=1, \ldots, n$,  $m\in\mathbb Z$, and central element $c$ with the product
$$[X_m, Y_n]=t^{m+n}[X,Y]+ \delta_{m+n,0}m(X|Y)c.$$

Let $\mathfrak a$ be a Lie algebra with a Cartan subalgebra $H$
and root system $\D$. Denote by $U(\mathfrak a)$ the universal
enveloping algebra of $\mathfrak a$.
 A closed subset $P\subset \D$ is called a
partition if $P\cap (-P)=\emptyset$ and $P\cup (-P)=\D$. If
$\mathfrak a$ is finite-dimensional then every partition
corresponds to a choice of positive roots in $\D$ and all
partitions are conjugate by the Weyl group. The situation is
different in the infinite-dimensional case. If $\mathfrak a$ is an
affine Lie algebra then partitions are divided into a finite
number of Weyl group orbits (cf. \cite{MR89m:17032}, \cite{MR1627814}).

Given a partition $P$ of $\D$ we define a {\em Borel} subalgebra
$\mathfrak b_P\subset \mathfrak a$ generated by $H$ and the root
spaces ${\mathfrak a}_{\alpha}$ with $\alpha \in P$. All Borel
subalgebras are conjugate in the finite-dimensional case. A {\em
parabolic} subalgebra is a subalgebra that contains a Borel
subalgebra. If $\mathfrak p$ is a parabolic subalgebra of a
finite-dimensional $\mathfrak a$ then $\mathfrak p=\mathfrak
p_0\oplus \mathfrak p_+$ where $\mathfrak p_0$ is a reductive Levi
factor and $\mathfrak p_+$ is a nilpotent subalgebra. Parabolic
subalgebras correspond to a choice of a basis $\pi$ of the root
system $\D$ and a subset $S\subset \pi$. A classification of all
Borel subalgebras in the affine case was obtained in \cite{MR1627814}. In
this case not all of them  are  conjugate but there exists a
finite number of conjugacy classes. These conjugacy classes are
parametrized by parabolic subalgebras of the underlined
finite-dimensional Lie algebra. Namely, let  $\mathfrak
p=\mathfrak p_0\oplus \mathfrak p_+$ a parabolic subalgebra of
$\G$ containing a fixed Borel subalgebra $\mathfrak b$ of
$\mathfrak g$. Define
$$
B_{\mathfrak p}=\left(\mathfrak p_+\otimes \C[t,t^{-1}]\right)\oplus \left(\mathfrak p_0\otimes t\C[t]\right)
\oplus \mathfrak b \oplus \C c\oplus \C d.
$$
For any Borel
subalgebra $\mathfrak B$ of $\hat{\mathfrak g}$ there exists a
parabolic subalgebra $\mathfrak p$ of $\mathfrak g$ such that
$\mathfrak B$ is conjugate to $B_{\mathfrak p}$.

When $\mathfrak p$ coincides with $\mathfrak g$, i.e. $\mathfrak
p_+=0$, the corresponding Borel subalgebra $B_{\mathfrak g}$ is
the {\em standard} Borel subalgebra defined by the choice of
positive roots in $\hat{\mathfrak g}$.  Another extreme case is
when $\mathfrak p_0=\H$. This corresponds to the {\em natural}
Borel subalgebra $B_{nat}$ of $\hat{\G}$ considered in \cite{MR89m:17032}.

Given a parabolic subalgebra $\mathfrak p$ of $\mathfrak g$ let
$\lambda: B_{\mathfrak p}\rightarrow \C$ be a $1$-dimensional
representation of $B_{\mathfrak p}$. Then one defines an induced
{\em Verma type} $\hat{\mathfrak g}$-module $$M_{\mathfrak
p}(\lambda)=U(\hat{\mathfrak g})\otimes_{U(B_{\mathfrak p})}\C.$$
The module $M_{\mathfrak g}(\lambda)$ is the classical Verma
module with highest weight $\lambda$ \cite{K}. In the case of
natural Borel subalgebra we obtain {\em imaginary} Verma modules
studied in \cite{MR95a:17030}. Note that the module $M_{\mathfrak
p}(\lambda)$ is $U(\mathfrak p_-)$-free, where $\mathfrak p_-$ is
the opposite subalgebra to $\mathfrak p_+$. The theory of Verma
type modules was developed in \cite{MR1627814}. It follows immediately
from the definition that, unless it is a classical Verma module, a
Verma type module with highest weight $\lambda$ has a unique
maximal submodule, it has both finite and infinite-dimensional
weight spaces and it can be obtained using the parabolic induction
from a standard Verma module $M$ with highest weight $\lambda$
over a certain affine Lie subalgebra. Moreover, if the central
element $c$ acts non-trivially on such Verma type module, then the
structure of this module is completely determined by the structure
of module $M$, which is well-known (\cite{MR1627814}, \cite{MR95d:17026}).

Let ${\mathfrak n}^{\pm}=\oplus_{\alpha\in {\Delta}^+}{\mathfrak
g}_{\pm \alpha}$. Denote
 ${\mathfrak
n}^{\pm}_r={\mathfrak n}^{\pm}\cap {\mathfrak g}_r$, ${\mathfrak
n}^{\pm}(r)={\mathfrak n}^{\pm}\setminus {\mathfrak n}^{\pm}_r$,
$$\bar{B}_r=L({\mathfrak n}^+(r))\oplus ({\mathfrak n}^{+}_r\otimes
\mathbb C[t]) \oplus (({\mathfrak n}^-_r)\oplus {\mathfrak
H})\otimes \mathbb C[t]t).$$
 Then
 $B_r=\bar{B}_r\oplus \hat{\mathfrak H}$ is a Borel subalgebra of
$\hat{\mathfrak g}$ for any $0\leq r\leq n$.

Given $\tilde{\lambda}\in \hat{\mathfrak H}^*$  the corresponding
Verma type module is $$M_r(\tilde{\lambda})=U(\hat{\mathfrak
g})\otimes_{U(B_r)}\mathbb C v_{\tilde{\lambda}}$$ where
$\bar{B}_r  v_{\tilde{\lambda}}=0$ and $h v_{\tilde{\lambda}}=
\tilde{\lambda}(h) v_{\tilde{\lambda}}$ for all $h\in
\hat{\mathfrak H}$.

When $r=n$ it gives us the usual Verma module construction. If $r=0$ we
get an imaginary Verma module.

Let $\tilde{\lambda}_r=\tilde{\lambda}|_{\hat{\mathfrak H}_r}$.
Verma type module  $M_r(\tilde{\lambda})$ contains a
$\hat{\mathfrak g}_r$-submodule
$M(\tilde{\lambda}_r)=U(\hat{\mathfrak g}_r)(1\otimes
v_{\tilde{\lambda}})$ which is isomorphic to a usual Verma module
for $\hat{\mathfrak g}_r$.

\begin{thm}[\cite{C}, \cite{FS}]\label{vermatype}
Let $\tilde{\lambda}(c)\neq 0$. Then the submodule structure of
$M_r(\tilde{\lambda})$ is completely determined by the submodule
structure of $M(\tilde{\lambda}_r)$. In particular,
$M_r(\tilde{\lambda})$ is irreducible if $M(\tilde{\lambda}_r)$ is
irreducible.

\end{thm}

Let $V$ be a weight $\mathfrak a$-module, i.e. $V=\oplus_{\mu\in
H^*}  V_{\mu}$, $V_{\mu}=\{v\in V|hv=\mu(h)v, \forall h\in H\}$.
Suppose that $\dim V_{\mu}\le \infty$ for all $\mu$. If $\mathfrak
a$ is a Kac-Moody Lie algebra (finite or affine) with Serre
generators $e_i$'s and $f_i$'s, then denote by $w$ an
anti-involution on $\mathfrak a$ which permutes $e_i$ with $f_i$
for all $i$, and which is the identity on $H$. Consider the $\mathfrak
a$-module
$$V^*=\oplus_{\mu\in
H^*}  V_{\mu}^*,$$ where $V_{\mu}^*$ is a dual subspace of
$V_{\mu}$ and the structure of $\mathfrak a$-module is given by:
$(xf)(v)=f(w(x)v)$. Such modules are called {\em contragradient}.
The advantage of considering these modules versus the whole dual
modules of $V$ is that $V$ and its contragradient module belong to
the Bernstein-Gelfand-Gelfand category $\mathcal O$
simultaneously.

\section{Geometric realizations}
\subsection{Finite-dimensional case}
Let  $\G$ be a finite-dimensional simple Lie algebra, $\mathfrak
b$ a Borel subalgebra, $\G=Lie G$, $\mathfrak b=Lie B$. Then the
group $G$ acts transitively on the {\em flag variety} $X=G/B$.
 Let
 $\mathcal D_X$ be the sheaf of differential operators on $X$
with regular coefficients. If $\mathcal M$ is a $\mathcal
D_X$-module, then the global sections $\Gamma(X, \mathcal M)$ have
a structure of a $\G$-module.

 If
$\lambda\in \H^*$, then the Verma module $M(\lambda)$ with the
highest weight $\lambda$ admits the central character which we
denote by $\xi_{\lambda}.$ Let $\xi: Z(\G)\rightarrow \C$ and
$\lambda\in \H^*$ be such that $\xi=\xi_{\lambda}$. Consider
$U_{\lambda}=U(\G)/U(\G)Ker \xi$. Then there is an isomorphism of
algebras
$$U_{\lambda}\simeq \Gamma(X, \mathcal D_{\lambda}),$$
where $\mathcal D_{\lambda}$ is a {\em twisted sheaf} of
differential operators  on $X$ introduced by Beilinson and
Bernstein (\cite{MR55:2941}).  Moreover, there is an equivalence
between the category of $U_{\lambda}$-modules and the category of
quasi-coherent $\mathcal D_{\lambda}$-modules on $X$
\cite{MR55:2941}.

 Let $\G=\mathfrak n_-\oplus \H\oplus \mathfrak n_+$, $\mathfrak b=\mathfrak n_-\oplus \H$,
 $\mathfrak n_{\pm}=\text{Lie} N_{\pm}$. Then $X=G/B$ has a decomposition into open Schubert cells:
$X=\cup_{w\in W}C(w)$, where $C(w)=B_+wB_-/B_-$, $W=N(T)/T$ is the
Weyl group and $T=B_+/N_+$. The subgroup $N_+$ acts on $X$, and
the largest orbit $\mathcal U$ of this action can be identified
with proper $N_+$. The Lie algebra $\G$ can be mapped into vector
fields on $X$ and hence on $\mathcal U$. Thus $\G$ can be embedded
into the differential operators on $\mathcal U$ of degree $\leq
1$.  Note that the ring of regular functions $\mathcal O_{\mathcal
U}$ on $\mathcal U$ is  a polynomial ring in $m=|\Delta_+|$
variables and hence $\G$ has an embedding into the Weyl algebra
$\mathcal A_m$ generated by $x_1, \ldots, x_n$ and partial
derivatives $\partial_1, \ldots,
\partial_n$. If $\xi:\Z(\G)\rightarrow \C$ is the central
character of $\G$ then the quotient $U(\G)/(\text{Ker} \xi)U(\G)$ can be
embedded into $A_n$, $n=(1/2)(\dim \G - \text{rank} \G)$ \cite{MR0374214},
providing a realization of $\G$ as differential operators acting on the Fock space $\C[x_1,\ldots,
x_n]$. A different approach was suggested by Khomenko \cite{Kh},
who showed that the quotient $U(\mathfrak{gl}(n))/(\text{Ker}\, \xi)U(\mathfrak{gl(n)})$ can be
embedded into a certain localization of $A_m$, $m=n(n+1)/2$, using
the theory of Gelfand-Tsetlin modules \cite{MR1308982}.

The embedding above of the Lie algebra $\G$ into $\mathcal A_m$
induces the structure of a $\G$-module on  $\mathcal O_{\mathcal
U}$. In fact, a $\G$-module $\mathcal O_{\mathcal U}$ is
isomorphic to a contragradient module $M^*(0)$ with trivial
highest weight.

For a general $\lambda\in \H^*$, $\Gamma(\mathcal U, \mathcal
D_{\lambda})\simeq M^*(\lambda)$ is a $\G$-module under the
inclusion
$$\G\subset \Gamma(G/B_-, \mathcal D_{\lambda})$$
 (Remark 10.2.7 in \cite{MR1849359}). In order to obtain a
geometric realization of Verma modules one needs to consider the
minimal $1$-point orbit of $N_+$ on $X$. Choosing another orbit of
$N_+$ gives a {\em twisted} Verma module parametrized by the
elements of the Weyl group. These modules have the same character
as corresponding Verma modules.

\subsection{Affine case}
Let $\hat{\G}$ be the non-twisted affine Lie algebra associated
with $\G$.
 Consider a Cartan decomposition $\G=\mathfrak n_-\oplus
\H\oplus \mathfrak n_+$ and  a Borel subalgebra $\mathfrak
b_{\pm}=\mathfrak n_{\pm}\oplus \H$. Denote
$$\hat{\mathfrak n}_{\pm}=(\mathfrak n_{\pm}\otimes 1)\oplus (\G\otimes t^{\pm}\C[t^{\pm}]),$$
$\hat{\mathfrak b}_{\pm}=\hat{\mathfrak n}_{\pm}\oplus \H\otimes
\C[t]$. Let $\hat{G}$, $\hat{N}_{\pm}$ and $\hat{B}_{\pm}$ be Lie
groups corresponding to $\hat{\G}$, $\hat{\mathfrak n}_{\pm}$ and
$\hat{\mathfrak b}_{\pm}$ respectively.
 Then a scheme of  infinite type $X=\hat{G}/\hat{B}_{-}$ which has a
 structure of a splits into $\hat{N}_{+}$-orbits of
 finite codimension, parametrized by the affine Weyl group. There
 is an analogue of a big cell $\hat{\mathcal U}$ in $X$ which is a
 projective limit of affine spaces, and hence, the ring of regular
 functions  $\mathcal O_{\hat{\mathcal U}}$ on $\hat{\mathcal U}$ is
 a polynomial ring in infinitely many variables. Thus $\hat{\G}$ acts
on it by differential operators
 providing a realization
 for the contragradient Verma module with zero highest weight.
Global sections of more general $\hat{N}_{+}$-equivariant sheaves
on $X$ will produce an
arbitrary highest weight. Other $\hat{N}_{+}$-orbits in $X$
 correspond to twisted contragradient Verma modules.
Standard Verma modules
 can be obtained by considering $\hat{N}_{+}$-orbits on
  $\hat{G}/{\hat{B}_{+}}$.

\subsection{First free field realization} In the previous section we
considered the case of classical Verma modules for affine Lie
algebras. Consider  now the natural Borel subalgebra  $\mathfrak
b_{nat}=\mathfrak n_-\otimes \C[t,t^{-1}]\oplus \H\otimes
\C[t^{-1}]$ of $\hat{\mathfrak g}$ and the corresponding Borel
subgroup $\mathfrak B_{nat}$. Let $X=\hat{G}/\mathfrak B_{nat}$.
which is {\em semi-infinite} manifold \cite{MR1849359},
\cite{MR1223226}, \cite{MR1713301}. We can consider the
$\hat{N}_+$-orbits on it and, in particular, $\hat{N_+}$ can be
viewed as an analogue of the big cell $\mathcal U$ in $G/B_-$.

Applying the same argument as in the previous section one can
obtain an embedding of $\hat{\G}$ into the Weyl algebra in
infinitely many variables and hence a realization of our algebra
in the Fock space $\C[x_n, n\in \Z]$. For example, consider the
case $\G=\mathfrak{sl}(2)$. Then
$$
e_n=e\otimes t^n, \, h_n=h\otimes t^n, \, f_n=f\otimes t^n, n\in \Z
$$
form a basis of $\G\otimes \C[t, t^{-1}]$. Then we have the
following embedding:
$$e_n\mapsto\frac{\partial}{ \partial x_n}, \,\,\, h_n\mapsto -2\sum_{m\in \Z}x_m\frac{\partial}{\partial x_{n+m}},\,\,\,
f_n\mapsto -\sum_{m,k\in \Z}x_mx_k\frac{\partial}{\partial x_{n+m+k}}.
$$ 
Note that
the differential operators corresponding to $f_n$ are not
well-defined on $\C[x_n, n\in \Z]$ (they take values in some
formal completion of $\C[x_n, n\in \Z]$). One way to deal with
this problem is to apply the anti-involutions:
$$
e_n\leftrightarrow f_n, \, \, \, \, h_n\leftrightarrow h_n; \,\,\,\,
 x_n\leftrightarrow\partial x_n, \,\,n\in
\Z
$$ 
which gives the following formulas:
$$
f_n\mapsto\frac{\partial}{\partial x_n}, \,\,\, h_n\mapsto -2\sum_{m\in \Z}x_{n+m}\frac{\partial}{\partial x_{m}},\,\,\,
e_n\mapsto -\sum_{m,k\in \Z}x_{n+m+k}\frac{\partial}{\partial x_{m}}\frac{\partial}{\partial
x_{k}}.
$$

These formulas define the {\em first free field realization} of
$\hat{sl}(2)$ in the polynomial ring $\C[x_m, m\in \Z]$.  This
module is, in fact, a quotient $M(0)$ of the imaginary Verma
module with trivial highest weight by a submodule generated by the
elements $h_n\otimes 1$, $n<0$. Similar formulas for an arbitrary
highest weight with a trivial action of the central element were
obtained by Jakobsen and Kac \cite{MR89m:17032} using analytic approach.
Namely, if $\mu$ denotes a finite measure on the circle $S^1$, not
concentrated in a finite number of points, and
$\lambda_m=\int_{S^1}z^m d\mu$ then
$$
f_n\mapsto\frac{\partial}{\partial x_n}, \,\,\,\,\, h_n\mapsto -\lambda_n -2\sum_{m\in \Z}x_{n+m}
\frac{\partial}{\partial x_{m}}
$$
$$
e_n\mapsto -\sum_{m,k\in \Z}x_{n+m+k}\frac{\partial}{ \partial x_{m}}\frac{\partial}{\partial
x_{k}}-\sum_{m\in \Z}\lambda_{n+m}\frac{\partial}{\partial x_m}
$$ 
gives a boson type
realization of $M(\lambda)$, where $\lambda(c)=0$ and
$\lambda(h_0)=-\lambda_0$. This module is irreducible if
$\lambda_0\neq 0$ \cite{MR95a:17030}.

To get a realization of the imaginary Verma module for
$\hat{\mathfrak{sl}}(2)$ with a non-trivial central action Bernard and Felder
used the Borel-Weil construction. Let $\hat{B}_-$ be the Borel
subgroup of the loop group $\hat{SL}(2)$ corresponding to a Borel
subalgebra $\mathfrak b_{nat}$. Then we presume $\hat{B}_-$
consists of the elements of the form
$$
\exp(\sum_{n\in \Z} x_n e_n)\exp(\sum_{m>0}y_mh_m),
$$
where $x_n, y_m$ are coordinate functions. Consider a one
dimensional representation  $\chi: \hat{B}_-\rightarrow \C$, where
$c$ acts by scalar $K$,  $h_0$ acts by scalar $J$ and all other
elements act trivially. Then the group $\hat{SL}(2)$ acts on the
sections of the line bundle ($\mathcal L_{\chi}, g)$ over
$\hat{SL}(2)/\hat{B}_-$ by
$$(g_1 f)(g_2)=f(g_1^{-1}g_2),$$
$g_i\in \hat{SL}(2)$, where
$$\mathcal L_{\chi}=\hat{SL}(2) \times_{\hat{B}_-} \C$$ and
$g:\hat{SL}(2) \times_{\hat{B}_-} \C \rightarrow
\hat{SL}(2)/\hat{B}_-$ such that $(x,z)\mapsto x \hat{B}_-$.

Differentiating this action to an action of the Lie algebra
$\hat{\mathfrak g}$ and applying two anti-involutions
$$
e_n\leftrightarrow -f_{-n}, \,\, h_n\leftrightarrow h_{-n}, \,\,\, c\leftrightarrow c
$$
and
$$
x_{-n}\leftrightarrow\partial/ \partial x_n, \, \, \, \, \, y_k\leftrightarrow -\partial/\partial y_k.
$$
we obtain the following boson realization of $\hat{\mathfrak g}$
in the Fock space $\C[x_m, m\in \Z]\otimes \C[y_n, n>0]$:
$$
f_n\mapsto  x_n, \,\,\, h_n\mapsto -2\sum_{m\in \Z}x_{m+n}\partial/\partial x_{m}
+\delta_{n<0}y_{-n} + \delta_{n>0}2nK\partial/\partial y_{n}
+\delta_{n,0}J,
$$
$$
e_n\mapsto -\sum_{m,k\in \Z}x_{k+m+n}\frac{\partial^2}{\partial x_k \partial x_{m}}
+\sum_{k>0}y_{k}\frac{\partial}{\partial x_{-k-n}}+ 2K\sum_{m>0}m\frac{\partial^2}{\partial
y_m\partial x_{m-n}}+(Kn+J)\frac{\partial}{\partial x_{-n}}.$$

 This module is irreducible if and only if $K\neq 0$. If we
 let $K=0$ and quotient out the submodule generated by $y_m, m>0$
then the factor module is irreducible if and only if $J\neq 0$ (cf.
\cite{MR95a:17030}). This construction has been generalized for all affine
Lie algebras in \cite{C2} providing a realization of imaginary
Verma modules.

\section{Boson type realizations}

The generators of Weyl algebras are called {\it free bosons}.
Hence any embedding of the affine Lie algebra $\hat{\G}$ into a
Weyl algebra leads to a boson type realization.

\subsection{Second free field realization} There is another way to
correct the formulas above using the construction of Wakimoto
modules \cite{W}.

Denote $a_n=\partial/\partial x_n$, $a_n^*=x_{-n}$ and consider formal
power series $$a(z)=\sum_{n\in \Z}a_n z^{-n-1}, \,\,\,\,
a^*(z)=\sum_{n\in \Z}a_n^* z^{-n}.$$ Series $a(z)$ and $a^*(z)$
are called {\em formal distributions}. It is easy to see that
$[a_n, a_m^*]=\delta_{n+m,0}$ and all other products are zero. The
formulas in the $sl(2)$ case can be rewritten as follows:
$$e(z)\mapsto a(z), \,\,\, h(z)\mapsto -2a^*(z)a(z), \,\,\,f(z)\mapsto -a^*(z)^2 a(z),$$
where $g(z)=\sum_{n\in \Z}g_n z^{-n-1}$ for $g\in \{e,f,h\}$. This
realization is not well-defined since the annihilation and
creation operators are in a wrong order.  It becomes well-defined
after the application of two anti-involutions described above.
Then the formulas read:
$$f
(z)\mapsto a(z), \,\,\, h(z)\mapsto 2a(z)a^*(z),
\,\,\,f(z)\mapsto -a(z) a^*(z)^2,
$$
where $a_n$ and $a_n^*$ have
the following meaning now $a_n=x_n$, $a_n^*=-\partial x_{-n}$.
This is our quotient of the imaginary Verma module.

A different approach was suggested by Wakimoto (\cite{W}) who used
the technique of {\em normal ordering}.  Denote
$$
a(z)_-=\sum_{n<0}a_nz^{-n-1}, \,\,\, a(z)_+=\sum_{n\geq
0}a_nz^{-n-1}
$$
and define the normal ordering as follows
$$:a(z)b(z):=a(z)_-b(z)+b(z)a_+(z).$$

Let now $$a_n=\begin{cases}
x_n,  & n<0\\
\frac{\partial}{\partial x_n},  & n\geq 0,
\end{cases}
\,\,\, a^*_n=\begin{cases}
x_{-n},  & n\leq0\\
-\frac{\partial}{\partial x_{-n}},  & n>0,
\end{cases}
\,\,\, b_m=\begin{cases}
m\frac{\partial}{\partial y_m},  & m\geq 0\\
y_{-m},  & m<0.
\end{cases}
$$
Here  $[a_n, a^*_m]=[b_n, b_m]=\delta_{n+m,0}$.

It was hown in \cite{W} that the formulas
$$
c\mapsto K, \,\,\, e(z)\mapsto a(z), \,\,\, h(z)\mapsto -2:a^*(z)a(z): +b(z),
$$
$$
f(z)\mapsto -:a^*(z)^2 a(z): +K\partial_za^*(z) +a^*(z)b(z)
$$
define the action of the affine $\hat{\mathfrak{sl}}(2)$  on the space
$\C[x_n, n\in \Z]\otimes \C[y_m, m>0]$. This boson type
realization is called the {\em second free field realization}
giving the celebrated {\em Wakimoto modules}. For an arbitrary
affine Lie algebra Wakimoto modules were constructed by Feigin and
Frenkel \cite{MR92f:17026}, \cite{MR92d:17025}. Generically they
are isomorphic to Verma modules. Wakimoto modules can be viewed as
infinite twistings of Verma modules.

\begin{rem}
We see that the semi-infinite variety $\hat{G}/\mathfrak B_{nat}$
gives rise to boson type realizations of Imaginary Verma modules
(the first free field realization) and of Wakimoto modules (the
second free field realization). In fact, one can construct the
whole family of "other" free field realizations \cite{MR2065632}.
\end{rem}


\subsection{ Oscillator algebras}

Let $\hat{\mathfrak a}$ be the infinite
dimensional Heisenberg algebra with generators
$a_{ij,m}$, $a_{ij,m}^*$, and $\mathbf 1$,  $1\leq i\leq j\leq n$  and
$m\in
\mathbb Z$, subject to the relations
\begin{align*}
[a_{ij,m},a_{kl,n}]&=[a_{ij,m}^*,a^*_{kl,n}]=0, \\
[a_{ij,m},a^*_{kl,n}]&=\delta_{ik}\delta_{jl}
\delta_{m+n,0}\mathbf 1, \\
[a_{ij,m},\mathbf 1]&=[a^*_{ij,m},\mathbf 1]=0.
\end{align*}
Such an algebra has a representation $\tilde{\rho}:\hat{\mathfrak
a}\to {\mathfrak{gl}}(\mathbb
C[\mathbf x])$ where
\begin{align*}   \mathbb C[\mathbf x]&:=
       \mathbb C[x_{ij,m}|i,j,m\in \mathbb Z,\,1\leq i\leq j\leq
    n]
\end{align*}
denotes the algebra over $\mathbb C$ generated by the indeterminates
$x_{ij,m}$ and $\tilde{\rho}$ is
defined by
\begin{align*}
 \tilde{\rho}( a_{ij,m}):&=\begin{cases}
  \partial/\partial
x_{ij,m}&\quad \text{if}\quad m\geq 0,\enspace\text{and}\enspace  j\leq r
\\ x_{ij,m} &\quad \text{otherwise},
\end{cases}
 \\
\tilde{\rho}(a_{ij,m}^*):&=
\begin{cases}x_{ij,-m} &\enspace \text{if}\quad m\leq
0,\enspace\text{and}\enspace j\leq r \\ -\partial/\partial
x_{ij,-m}&\enspace \text{otherwise}. \end{cases}
\end{align*}
and $\tilde{\rho}(\mathbf 1)=1$.
In this case
$\mathbb C[\mathbf x]$ is an
$\hat{\mathfrak a}$-module generated by $1=:|0\rangle$, where
$$
a_{ij,m}|0\rangle=0,\quad m\geq  0 \enspace\text{and}\enspace j\leq r,
\quad a_{ij,m}^*|0\rangle=0,\quad m>0\enspace\text{or}\enspace j>r.
$$
Let $\hat{\mathfrak a}_r$ denote the subalgebra generated by $a_{ij,m}$
and
$a_{ij,m}^*$ and $\mathbf 1$, where $1\leq i\leq j\leq r$ and
$m\in\mathbb
Z$. If $r=0$, we set $\hat{\mathfrak a}_r=0$.

Let $A_n=((\alpha_i|\alpha_j))$ be the Cartan matrix for $\mathfrak{sl}(n+1,\mathbb C)$
and
let $\mathfrak B$ be the matrix whose entries are
$$
\mathfrak B_{ij}:=(\alpha_i|\alpha_j)(\gamma^2
    -\delta_{i>r}\delta_{j>r}(r+1)
    +\frac{r}{2}\delta_{i,r+1}\delta_{j,r+1})
$$

where $$
\delta_{i>r}= \begin{cases}
1& \quad \text{if}\quad i>r,\\
0& \quad \text{otherwise}.
\end{cases}$$

In other words
$$
\mathfrak B:=\gamma^2A_n -(r+1)
\begin{pmatrix} 0  & 0 \\ 0 & A_{n-r} \end{pmatrix}+rE_{r+1,r+1}.
$$
Using cofactor expansion along the $r+1$-st row one can show
$$
\det \mathfrak B=(n+1)\gamma^{2r}(\gamma^2-r-1)^{n-r}
$$
where we define $\gamma^{2r}=1$ when $r=0$ and $\gamma=0$.
Thus $\mathfrak B$ is degenerate if and only if either $\gamma=0$ or $\gamma^2=r+1$.
We also have the Heisenberg Lie algebra $\hat{\mathfrak b}$ with generators $b_{im}$,
$1\leq
i\leq n$, $m\in\mathbb Z$, $\mathbf 1$, and relations
$[b_{im},b_{jp}]=m\,\mathfrak B_{ij}\delta_{m+p,0}\mathbf 1$ and
$[b_{im},\mathbf 1]=0$.

 For each
$1\leq i\leq n$ fix $\lambda_i\in\mathbb C$ and let
$\lambda=(\lambda_1, \ldots, \lambda_n)$.
 Then the algebra  $\hat{\mathfrak b}$
has
a representation
$\rho_{\lambda}:\hat{\mathfrak b}\to\End(\mathbb C[\mathbf y]_\lambda)$ where
\begin{align*}   \mathbb C[\mathbf y]&:=
       \mathbb C[y_{i,m}|i,m\in \mathbb N^*,\,1\leq i\leq
    n]
\end{align*}
  and $\rho_{\lambda}$ is
defined
on
$\mathbb C[\mathbf y]$ defined by
$$
\rho_{\lambda}(b_{i0})= \lambda_i, \quad \rho_{\lambda}(b_{i,-m})=
     \mathbf e_i\cdot \mathbf y_{m},\quad
     \rho_{\lambda}(b_{im})=m\mathbf e_i \cdot  \,\frac{\partial}{\partial
\mathbf y_{m}}\quad
\text{for}\quad m>0
$$
and $\rho_{\lambda}(\mathbf 1)=1$.  Here
$$
\mathbf y_m=(y_{1m},\cdots, y_{nm}),\quad \frac{\partial}{\partial
\mathbf y_{m}}=\left(\frac{\partial}{\partial   y_{1m}}, \cdots,
\frac{\partial}{\partial y_{nm}} \right)
$$
and $\mathbf e_i$ are vectors in $\mathbb C^n$ such that $\mathbf
e_i\cdot
\mathbf e_j=\mathfrak B_{ij}$ where $\cdot $ means the usual dot
product.

Note that since $\mathfrak B_{ij}$ is symmetric, it is
orthogonally diagonalizable, (i.e. there exists an orthogonal matrix $P$
such that $P^t\mathfrak BP$ is a diagonal matrix) and hence we can find
vectors $\mathbf e_i$  in
$\mathbb C^n$ such that
$\mathbf e_i\cdot
\mathbf e_j=\mathfrak B_{ij}$. In fact for
$m>0$ and
$n<0$ we get
\begin{align*}
[b_{im},b_{jn}]& =m\delta_{m+n,0}\mathfrak
B_{ij}.
\end{align*}
(See also \cite{MR92d:17025}.)

We also define for $0\leq r\leq n$
\[
  \mathbb C_r[\mathbf y]:=
       \mathbb C[\mathbf e_i\cdot \mathbf y_{m}|i,m\in \mathbb N^*,\,1\leq i\leq
    r].
\]
If $\mathfrak B$ is non-degenerate, then $  \mathbb C_n[\mathbf
y]=  \mathbb C[\mathbf y]$. On the other hand if $\mathfrak B$ is
degenerate, then there is some choice to be made for the $\mathbf
e_i$.  For example if $\gamma=0$, $n=2$ and $r=1$, then
$$
\mathfrak B=\begin{pmatrix} 0 & 0 \\ 0 & -3 \end{pmatrix}
$$
and one can choose for example $\mathbf e_1=(0,0)$ and $\mathbf
e_2=(0,\sqrt{-3})$ or $\mathbf e_1=(0,0)$ and $\mathbf
e_2=(\sqrt{-3},0)$. In either case $\C_n[\mathbf y]\neq \C[\mathbf
y]$.

\subsection{Formal Distributions}
We need some more notation that will simplify some of the arguments
later.
This notation follows roughly \cite{MR99f:17033} and \cite
{MR2000k:17036}:  A {\it formal distribution} is an expression of the form
$$
a(z,w,\dots)=\sum_{m,n,\dots\in\mathbb Z}a_{m,n,\dots}z^mw^n
$$
where the $a_{m,n,\dots}$ lie in some fixed vector space $V$ and
$z,w,\dots $ are formal variables. We define $\partial a(z)=\partial_z
a(z)=\sum_nna_nz^{n-1}$. We also have expansion about zero:  there are
two canonical embeddings of vector spaces $\iota_{z,w}:\mathbb
C(z-w)\to
\mathbb C[[z,w]]$ and  $\iota_{w,z}:\mathbb C(z-w)\to
\mathbb C[[z,w]]$
where $
\iota_{z,w}(a(z,w))$ is formal Laurent series expansion in $z^{-1}$
and
$-\iota_{w,z}(a(z,w))$  is formal Laurent series expansion in $z$.
The {\it formal delta function}
$\delta(z-w)$ is the formal distribution
$$
\delta(z-w)=z^{-1}\sum_{n\in\mathbb Z}\left(\frac{z}{w}\right)^n
=\iota_{z,w}\left(\frac{1}{z-w}\right)-\iota_{w,z}\left(
\frac{1}{z-w}\right).
$$
For any sequence of elements $\{a_{m}\}_{m\in
\mathbb Z}$ in the ring $\End (V)$, $V$ a vector space,  the
formal distribution
\begin{align*}
a(z):&=\sum_{m\in\mathbb Z}a_{m}z^{-m-1}
\end{align*}
is called a {\it field}, if for any $v\in V$, $a_{m}v=0$ for $m\gg0$.
If $a(z)$ is a field, then we set
\begin{align*}
    a(z)_-:&=\sum_{m\geq 0}a_{m}z^{-m-1},\quad\text{and}\quad
   a(z)_+:=\sum_{m<0}a_{m}z^{-m-1}.
\end{align*}

 In particular
$$
\delta(z-w)_-=\iota_{z,w}\left(\frac{1}{z-w}\right),\quad
\delta(z-w)_+=-\iota_{w,z}\left(\frac{1}{z-w}\right).
$$
Note that
$$
-\partial_z\delta(z-w)=\partial_w\delta(z-
w)=\iota_{z,w}\left(\frac{1}{(z-w)^2}\right)
-\iota_{w,z}\left(\frac{1}{(z-w)^2}\right).
$$

 The {\it normal ordered product} of two distributions
$a(z)$ and
$b(w)$ (and their coefficients) is
defined by
\begin{equation}\label{normalorder}
\sum_{m\in\mathbb Z}\sum_{n\in\mathbb
Z}:a_mb_n:z^{-m-1}w^{-n-1}=:a(z)b(w):=a(z)_+b(w)+b(w)a(z)_-.
\end{equation}

For any $1\leq i\leq j\leq n$, we define
$$
a^*_{ij}(z)=\sum_{n\in\mathbb Z}a^*_{ij,n}z^{-n},\quad
a_{ij}(z)=\sum_{n\in\mathbb Z}a_{ij,n}z^{-n-1}
$$
and
$$
b_i(z)=\sum_{n\in\mathbb Z}b_{in}z^{-n-1}.
$$
In this case
\begin{align*}
[b_i(z),b_j(w)]
&=\mathfrak{B}_{ij}\partial_w\delta(z-w), \\
[a_{ij}(z),a^*_{kl}(w)]&
    =\delta_{ik}\delta_{jl}\mathbf 1\delta(z-w).
\end{align*}

Observe that $a_{ij}(z)$ for $j>r$ is not a field whereas $a_{ij}^*(z)$
is always a field.   We will call $a_{ij}(z)$ (resp. $a^*_{ij}(z)$)
a {\it pure creation} (resp. {\it annihilation}) {\it operator} if
$j>r$. Set
\begin{align*}
a_{ij}(z)_+&=a_{ij}(z),\quad a_{ij}(z)_-=0 \\
a_{ij}^*(z)_+&=0,\quad a_{ij}^*(z)_-=a_{ij}^*(z) ,
\end{align*}
if $j>r$.

Now we should point out that while $:a^1(z_1)\cdots a^m(z_m):$
is always defined as a formal series, we will only define $:a(z)
b(z)::=\lim_{w\to z}:a(z)b(w):$
for certain pairs
$(a(z),b(w))$.  For example
$$
:a_{ij}(z)a^*_{kl}(z):=\sum_{m\in\mathbb Z}\left(\sum_{n\in\mathbb
Z}:a_{ij,n}a^*_{kl,m-n}:\right)z^{-m-1}
$$
is well defined as an element in $\End(\mathbb C[\mathbf x]\otimes
\mathbb C[\mathbf y])[[z,z^{-1}]]$ for all $l>r$ (as
$\tilde{\rho}(a_{kl,m}^*):=-\partial/\partial x_{kl,-m}$ for $l>r$ )
or if both $l\leq r$ and
$j\leq r$ (see also the remarks after \thmref{realization}).

Then one defines recursively
\[
:a^1(z_1)\cdots a^k(z_k):=:a^1(z_1)\left(:a^2(z_2)\left(:\cdots
:a^{k-1}(z_{k-1}) a^k(z_k):\right)\cdots
:\right):,
\]
while normal ordered product
\[
:a^1(z)\cdots
a^k(z):=\lim_{z_1,z_2,\cdots, z_k\to
z} :a^1(z_1)\left(:a^2(z_2)\left(:\cdots :a^{k-1}(z_{k-1})
a^k(z_k):\right)\cdots
\right):
\]
will only be defined for certain $k$-tuples $(a^1,\dots,a^k)$.

Let
\begin{equation}\label{contraction}
\lfloor
ab\rfloor=a(z)b(w)-:a(z)b(w):= [a(z)_-,b(w)],
\end{equation}
(half of
$[a(z),b(w)]$) denote the {\it contraction} of any two formal distributions
$a(z)$ and $b(w)$ where $a(z)$, $b(z)$ are free fields or pure
creation or annihilation operators. For example if
$j,l\leq r$, then
\begin{align}
\lfloor a_{ij} a_{kl}^*\rfloor
&
=\sum_{m\geq 0}\delta_{ik}\delta_{jl}z^{-m-1}w^{m}
=\delta_{i,k}\delta_{j,l}\delta_-(z-w)
=\delta_{ik}\delta_{jl}
\,\iota_{z,w}\left(\frac{1}{z-w}\right)\\
\lfloor a_{kl}^*a_{ij}\rfloor
&
=-\sum_{n<0}\delta_{ik}\delta_{jl}z^{n}w^{-n-1}
=-\delta_{ik}\delta_{jl}\delta_+(w-z)=
\delta_{ik}\delta_{jl}\,\iota_{z,w}\left(\frac{1}{w-z}
\right).
\end{align}
 If $l>r$, then
\begin{align}
\lfloor a_{ij} a_{kl}^*\rfloor
&=[a_{ij}(z)_-,a_{kl}^*(w)]=0 \\
\lfloor a_{kl}^*a_{ij}\rfloor
&
=[a_{kl}^*(z)_-,a_{ij}(w)]
=-\delta_{ik}\delta_{jl}\delta(w-z).
\end{align}

\begin{thm}[Wick's Theorem, \cite{MR85g:81096}, \cite{MR99m:81001},
\cite{MR99f:17033}, or \cite{MR2065632}]  Let  $a^i(z)$ and $b^j(z)$ be formal
distributions with coefficients in the associative algebra
 $\End(\mathbb C[\mathbf x]\otimes \mathbb C[\mathbf y])$,
 satisfying
\begin{enumerate}
\item $[ \lfloor a^i(z)b^j(w)\rfloor ,c^k(x)_\pm]=[ \lfloor
a^ib^j\rfloor ,c^k(x)_\pm]=0$, for all $i,j,k$ and
$c^k(x)=a^k(z)$ or
$c^k(x)=b^k(w)$.
\item $[a^i(z)_\pm,b^j(w)_\pm]=0$ for all $i$ and $j$.
\item The products
$$
\lfloor a^{i_1}b^{j_1}\rfloor\cdots
\lfloor a^{i_s}b^{i_s}\rfloor:a^1(z)\cdots a^M(z)b^1(w)\cdots
b^N(w):_{(i_1,\dots,i_s;j_1,\dots,j_s)}
$$
have coefficients in
$\End(\mathbb C[\mathbf x]\otimes \mathbb C[\mathbf y])$ for all subsets
$\{i_1,\dots,i_s\}\subset \{1,\dots, M\}$, $\{j_1,\dots,j_s\}\subset
\{1,\cdots N\}$. Here the subscript
${(i_1,\dots,i_s;j_1,\dots,j_s)}$ means that those factors $a^i(z)$,
$b^j(w)$ with indices
$i\in\{i_1,\dots,i_s\}$, $j\in\{j_1,\dots,j_s\}$ are to be omitted from
the product
$:a^1\cdots a^Mb^1\cdots b^N:$ and when $s=0$ we do not omit
any factors.
\end{enumerate}
Then
\begin{align*}
:&a^1(z)\cdots a^M(z)::b^1(w)\cdots
b^N(w):= \\
  &\sum_{s=0}^{\min(M,N)}\sum_{i_1<\cdots<i_s,\atop
j_1\neq \cdots \neq j_s}\lfloor a^{i_1}b^{j_1}\rfloor\cdots
\lfloor a^{i_s}b^{i_s}\rfloor
:a^1(z)\cdots a^M(z)b^1(w)\cdots
b^N(w):_{(i_1,\dots,i_s;j_1,\dots,j_s)}.
\end{align*}
\end{thm}

We will also need the following two
results.
\begin{thm}[Taylor's Theorem, \cite{MR99f:17033}, 2.4.3]
\label{Taylorsthm}  Let
$a(z)$ be a formal distribution.  Then in the region $|z-w|<|w|$,
\begin{equation}
a(z)=\sum_{j=0}^\infty \partial_w^{(j)}a(w)(z-w)^j.
\end{equation}
\end{thm}

\begin{thm}[\cite{MR99f:17033}, Theorem 2.3.2]\label{kac} Let $a(z)$ and $b(z)$
be formal distributions with coefficients in the associative algebra
 $\End(\mathbb C[\mathbf x]\otimes \mathbb C[\mathbf y])$.   The
following are equivalent
\begin{enumerate}[(i)]
\item
$\displaystyle{[a(z),b(w)]=\sum_{j=0}^{N-1}\partial_w^{(j)}
\delta(z-w)c^j(w)}$, where $c^j(w)\in \End(\mathbb C[\mathbf x]\otimes \mathbb
C[\mathbf y])[[w,w^{-1}]]$.
\item
$\displaystyle{\lfloor
ab\rfloor=\sum_{j=0}^{N-1}\iota_{z,w}\left(\frac{1}{(z-w)^{j+1}}
\right)
c^j(w)}$.
\end{enumerate}\label{Kacsthm}
\end{thm}

In other words the singular part of the {\it operator product
expansion}
$$
\lfloor
ab\rfloor=\sum_{j=0}^{N-1}\iota_{z,w}\left(\frac{1}{(z-w)^{j+1}}
\right)c^j(w)
$$
completely determines the bracket of mutually local formal
distributions $a(z)$ and $b(w)$.   One writes
$$
a(z)b(w)\sim \sum_{j=0}^{N-1}\frac{c^j(w)}{(z-w)^{j+1}}.
$$
For example
$$
b_i(z)b_j(w)\sim \frac{\delta_{ij}}{(z-w)^2}.
$$

\subsection{Intermediate Wakimoto modules}

Define
$$
E_{i}(z)=\sum_{n\in\mathbb Z}E_{in}z^{-n-1},\enspace
F_{i}(z)=\sum_{n\in\mathbb Z}F_{in}z^{-n-1}, \enspace
H_i(z)=\sum_{n\in\mathbb
Z} H_{in}z^{-n-1},\quad 1\leq i\leq n.
$$
The defining relations between the generators
of $\hat{\mathfrak g}$ can be written as follows

\begin{align*}
[H_i(z),H_j(w)]&=(\alpha_i|\alpha_j)c\partial_w \delta(w-z)
\tag{R1}\\
   [H_i(z),E_j(w)]&=(\alpha_i|\alpha_j)E_j(z)\delta(w-z)\tag{R2} \\
   [H_i(z),F_j(w)]&=-(\alpha_i|\alpha_j)F_j(z)\delta(w-z)\tag{R3} \\
   [E_i(z),F_j(w)]&=\delta_{i,j}(H_i(z)\delta(w-z)+c\partial_w
     \delta(w-z))\tag{R4} \\
   [F_i(z),F_j(w)]&=[E_i(z),E_j(w)]=0\quad\text{if}\quad
(\alpha_i|\alpha_j)\neq -1\tag{R5} \\
   [F_i(z_1),F_i(z_2),F_j(w)]&=[E_i(z_1),E_i(z_2),E_j(w)]=0
\quad\text{if}\quad (\alpha_i|\alpha_j)=
   -1\tag{R6}
   \end{align*}
where $[X,Y,Z]:=[X,[Y,Z]]$ is the Engel bracket for any three
operators
$X,Y,Z$.

Recall that $\mathbb C[\mathbf x]$ is an $\hat{\mathfrak
a}$-module with respect to the representation $\tilde{\rho}$ and
$\mathbb C[\mathbf y]$ is a $\hat{\mathfrak b}$-module with
respect to
 ${\rho}_{\lambda}$. In  \cite{MR2065632} we define a representation
$$\rho:\hat{\mathfrak g}\rightarrow
\mathfrak{gl}(\mathbb C[\mathbf x]\otimes \mathbb C[\mathbf y]).
$$
where we use the notation $\rho(X_m):=\rho(X)_m$, for $X\in\mathfrak
g$.
This is described in the following result:
\begin{thm}\label{realization}
Let $\lambda\in \mathfrak H^*$ and set $\lambda_i=\lambda(H_i)$.  The
generating functions

\begin{align}
\rho(F_{i})(z)&=a_{ii}+\sum_{j=i+1}^{n}a_{ij}a_{i+1,j}^*, \label{fs} \\ \notag \\
\rho(H_i)(z)&=2:a_{ii}a_{ii}^*:+\sum_{j=1}^{i-1}\left(:
         a_{ji}a_{ji}^*:-:a_{j,i-1} a_{j,i-1}^*:\right)  \\
       &\quad+\sum_{j=i+1}^{n}\left(:a_{ij}a_{ij}^*:
        - :a_{i+1,j}a_{i+1,j}^*:\right)
         + b_i, \notag \\  \notag   \\
\rho(E_i)(z)&= :a_{ii}^*\left(
       \sum_{k=1}^{i-1} a_{k,i-1}a_{k,i-1}^*-\sum_{k=1}^{i}
       a_{ki}a^*_{ki}\right): +\sum_{k=i+1}^n a_{i+1,k}a_{ik}^*
         -\sum_{k=1}^{i-1}a_{k,i-1}a_{ki}^* \\ \notag
&\quad-a_{ii}^*b_{i}-
     \left(\delta_{i> r}(r+1)+\delta_{i\leq
r}(i+1)-\gamma^2\right)\partial a_{ii}^*,\\
\notag \\
\rho(c)&=\gamma^2-(r+1)
\end{align}
define an action of the generators
$E_{im}$, $F_{im}$,
$H_{im}$, $i=1, \ldots, n$, $m\in \mathbb Z$ and $c$, on the Fock space
$\mathbb C[\mathbf x]\otimes
\mathbb C[\mathbf y]$.  In the above $a_{ij}$, $a_{ij}^*$ and $b_i$ denotes
$a_{ij}(z)$,
$a_{ij}^*(z)$ and $b_i(z)$ respectively.
\end{thm}

 \thmref{realization} defines a boson type realization of
$\hat{\mathfrak{sl}}(n+1, \mathbb C)$ and a module structure on
the Fock space $\mathbb C[\mathbf x]\otimes \mathbb C[\mathbf y]$
that depends on the parameter $0\leq r\leq n$. We called such
modules  {\it intermediate Wakimoto module} in \cite{MR2065632}. One can
easily see that \thmref{realization} defines also a boson type
realization of $\hat{\mathfrak{sl}}(n+1, \mathbb C)$ on the Fock
space $\mathbb C[\mathbf x]\otimes \mathbb C_n[\mathbf y]$ which
is different from the one above if $\mathfrak B$ is
non-degenerate. It is more convenient to work with such
realization, and we will call this module structure on $\mathbb
C[\mathbf x]\otimes \mathbb C_n[\mathbf y]$, the {\it intermediate
Wakimoto module} and denote it by $W_{n,r}(\lambda, \gamma)$.

 The intermediate Wakimoto
modules $W_{n,r}(\lambda, \gamma)$ have the property that the subalgebra
$\bar{B}_{r}$ annihilates the vector $1\otimes 1\in \mathbb C[\mathbf x]\otimes \mathbb C[\mathbf y]$, $h(1\otimes 1)=\lambda(h)(1\otimes 1)$ for all $h\in \mathfrak H$ and
$c(1\otimes 1)=(\gamma^2 - (r+1))(1\otimes 1)$.
Consider the $\hat{\mathfrak g}_r$-submodule $W=U(\hat{\mathfrak
g}_r)(1\otimes 1)
\simeq W_{r,r}(\lambda, \gamma)$ of $W_{n,r}(\lambda, \gamma)$.

\begin{rem} If $\lambda$ is generic then $W$ is isomorphic to the
Wakimoto module $W_{\lambda(r),\tilde{\gamma}}$
(\cite{MR92d:17025}) where $\lambda(r)=\lambda|_{\mathfrak H_r}$,
$\tilde{\gamma}=\gamma^2-(r+1)$.
\end{rem}

\begin{example}
If $\mathfrak B$ is non-degenerate then $W\simeq
W_{\lambda(r),\tilde{\gamma}}$. We also have an isomorphism in the
case $n=2$, $r=1$, $e_1=(0,0)$, $e_2=(\sqrt{-3}, 0)$. But there is
no isomorphism if we choose $e_2=(0, \sqrt{-3})$.

\end{example}

In particular we have the following corollary. Consider
$\tilde{\lambda}\in \hat{\mathfrak H}^*$ such that
$\tilde{\lambda}|_{\mathfrak H}=\lambda$,
$\tilde{\lambda}(c)=\gamma^2-(r+1)$, a Verma type module
$M_r(\tilde{\lambda})$ and its $\hat{\mathfrak g}_r$-submodule
$M(\tilde{\lambda}_r)$.

\begin{cor}\label{cor20}
Suppose that $M(\tilde{\lambda}_r)$ is irreducible and assume that
$\lambda(c)\neq 0$. Let $\tilde W=U(\hat{\mathfrak g})W$. Then
$M_r(\tilde{\lambda})$ isomorphic to $\tilde W$.
\end{cor}

\begin{proof}
Since $M(\tilde{\lambda}_r)$ is irreducible and $\lambda(c)\neq 0$
we conclude that $\mathfrak B$ is non-degenerate. Hence $W\simeq
W_{\lambda(r),\tilde{\gamma}}$ is essentially the Wakimoto module
for $\hat{\mathfrak g}_r$. Moreover, the Wakimoto module
$W_{\lambda(r),\tilde{\gamma}}$ is isomorphic to
$M(\tilde{\lambda}_r)$ in this case.  Indeed, the Verma module
$M(\tilde\lambda_r) $ is irreducible, so that the canonical map
$M(\tilde\lambda_r)\to W$, given by $1\mapsto 1\otimes 1$, is
injective. As $M(\tilde\lambda_r)$ and $W$ have the same character
formulae, this canonical map provides an isomorphism. We also have
a canonical map $\phi:M_{r}(\tilde\lambda_r)\to \tilde{W}$. This
map restricts to the canonical map $M(\tilde\lambda_r)\to W$, $W$
is contained in the image of $\phi$, hence $\phi$ is surjective.
Since $\tilde\lambda$ is generic, $M_{r}(\tilde\lambda)$ is
irreducible by Theorem \ref{vermatype} and thus $\phi$ must be an
isomorphism. Therefore $M_r(\tilde{\lambda})\simeq \tilde W$.
\end{proof}

It follows that Theorem \ref{realization} provides a boson type realization for generic Verma type modules.

\section{Generating Intermediate  Wakimoto modules}
Let $\hat{\mathfrak a}=\hat{\mathfrak{sl}}(r+1)$ with $\hat
{\mathfrak{a}} \subset\hat
{\mathfrak{g}}:=\hat{\mathfrak{sl}}(n+1)$. As above
$M_{r}(\tilde\lambda)$ will denote a Verma type module with
highest weight $\tilde\lambda$ and set
$$
W_{\mathfrak g}(\tilde\lambda):=\mathbb C[\mathbf x]\otimes \mathbb C_n[\mathbf y]
$$
with the action  defined by \thmref{realization} above.  When
$\mathfrak B$ is non-degenerate, $W_{\mathfrak g}(\tilde\lambda)$
is the intermediate Wakimoto module $W_{n,r}(\lambda, \gamma)$.
Sitting inside $W_{\mathfrak g}(\tilde\lambda)$ is a copy of
$W_{\mathfrak a}(\tilde{ \lambda}_r)$: For $1\leq k\leq n$, set
$\mathbb C_k[\mathbf x]:=\mathbb C[x_{ij,m}\,|\, 1\leq i\leq j\leq
k,m\in\mathbb Z]$ and
$$
W_{\mathfrak a}(\tilde{ \lambda}_r):=\mathbb C_r[\mathbf x]\otimes \mathbb C_r[\mathbf y].
$$
\begin{rem}
$W_{\mathfrak a}(\tilde{ \lambda}_r)\simeq W$ for a generic
$\tilde{\lambda}$.
\end{rem}

We will show that in fact $ W_{\mathfrak g}(\tilde\lambda)$ is
generated by $ W_{\mathfrak a}(\tilde{ \lambda}_r)$. Namely we
have

\begin{thm}\label{thm40}   $W_{\mathfrak g}(\tilde\lambda)=U(\hat{\mathfrak g})W_{\mathfrak a}(\tilde{ \lambda}_{r})$.
\end{thm}

\begin{proof}
First of all note that $L({\mathfrak n}^{-}(r))$ is generated by
$F_{i,m}$ with $r<i\leq n$ and $W_{\mathfrak a}(\tilde\lambda_r)$
is generated by $x_{ij,m}$ and $\mathbf e_j\cdot \mathbf y_{p}$
with $1\leq i\leq j\leq r$, $m\in\mathbb Z$ and $p\in\mathbb N$.

For $1\leq k\leq n$, set  $W_k:=\mathbb C_k[\mathbf x]\otimes
\mathbb C_n[\mathbf y]$. Let's first see that $W_r \subset
U(\widehat{\mathfrak g})W_{\mathfrak a}(\tilde\lambda_r)$.  Now an
arbitrary element in $W_r$ has the form
\[
\sum_j u_j\otimes v_j,\quad u_j\in \mathbb C_r[\mathbf x],
\quad v_j \in\mathbb C_n[\mathbf y]
\]
and so it suffices to show that any element of the form $u \otimes
v$, with $u\in \mathbb C_r[\mathbf x]$ and $v\in \mathbb
C_r[\mathbf y]$ monomials, is in $U(\widehat{\mathfrak
g})W_{\mathfrak a}(\tilde\lambda_r)$.  Write $v=v^rv_r$ with
$v_r\in  \mathbb C[\mathbf e_j\cdot \mathbf y_{m}\,|\, 1\leq j\leq
r,m\in\mathbb N]$ and $v^r\in \mathbb C[\mathbf e_j\cdot \mathbf
y_{m}\,|\, r< j\leq n,m\in\mathbb N]$ where we may assume
$v^r=(\mathbf e_{r}\cdot \mathbf y_{i_1})^{p_{i_1}}\cdots (\mathbf
e_{n}\cdot \mathbf y_{i_s})^{p_{i_s}}$ is a nonconstant monomial
(here $i_1,\dots ,i_s\in\mathbb N$). Recall
\begin{align*}
\rho(H_i)(z)&=2:a_{ii}a_{ii}^*:+\sum_{j=1}^{i-1}\left(:
         a_{ji}a_{ji}^*:-:a_{j,i-1} a_{j,i-1}^*:\right) \\
       &\quad+\sum_{j=i+1}^{n}\left(:a_{ij}a_{ij}^*:
        - :a_{i+1,j}a_{i+1,j}^*:\right)
         + b_i.
\end{align*}
For $i>r+1$ and $m\in\mathbb N$ those summands above  having
factors with $a^*_{ik,m}$ $k\geq r+1$ act as zero on $W_{\mathfrak
a}(\tilde\lambda_r)$ as these factors act by
$-\partial/\partial_{x_{ik,-m}}$.  Thus we get when restricted to
$W_{\mathfrak a}(\tilde\lambda_r)$ that $\rho(H_{im})=b_{im}$ and
this is given by left multiplication by $\mathbf e_i\cdot \mathbf
y_{m}$.  For $i=r+1$ and $m\in\mathbb N$, we have when restricted
to  $W_{\mathfrak a}(\tilde\lambda_r)$,
\begin{align*}
\rho(H_{r+1,m})&=-\sum_{j=1}^{r}\sum_{p\in\mathbb Z}:a_{jr,p} a_{jr,m-p}^*:
         + b_{r+1,m}.
\end{align*}
Altogether for $i\geq r+1$ and $m\in\mathbb N$, we can write
\begin{align*}
\rho(H_{im})&=X_{im}
         + b_{im}.
\end{align*}
where the first summand maps $ \mathbb C_r[\mathbf x]$ back into
itself and the second summand is left multiplication by $\mathbf
e_i\cdot \mathbf y_{m}$ ($X_{im}$ would be defined to be zero for
$i>r+1$).    As a consequence we get
\begin{align*}
u\otimes v&=u\otimes (\mathbf e_{r}\cdot \mathbf y_{i_1})^{p_{i_1}}\cdots  (\mathbf e_{n}\cdot \mathbf y_{i_s})^{p_{i_s}}v_r \\
&=b_{r+1,i_1}^{p_{i_1}}\cdots  b_{n, i_s}^{p_{i_s}}(u\otimes v_r )\\
&=(\rho(H_{r+1,i_1})-X_{r,i_1})^{p_{i_1}}
\cdots \left(\rho(H_{n, i_s})-X_{n,i_s}\right) ^{p_{i_s}}(u\otimes v_r ).
\end{align*}
Hence if we let $\mathcal H$ denote the abelian Lie algebra
generated by $H_{im}$ with $1\leq i\leq n$, $m\in\mathbb N$, then
\begin{equation}\label{ys}
u\otimes v\in U(\mathcal H)W_{\mathfrak a}(\tilde\lambda_r).
\end{equation}
and thus $W_r\subset U(\widehat{\mathfrak g})W_{\mathfrak a}(\tilde\lambda_r)$.

Fix $r\leq k<n$. Consider $u\otimes v\in W_{\mathfrak g}(\tilde\lambda)$ with
$$
u=x_{k+1,k+1,m_1}^{p} u_k,\quad u_k \in \mathbb C_k[\mathbf x],\quad m_1\in \mathbb Z,\quad p\in\mathbb N
$$  and $v\in \mathbb C_n[\mathbf y]$ and assume the induction hypothesis that $W_{k}\subset U(\widehat{\mathfrak g})W_{\mathfrak a}(\tilde\lambda_r)$.  (This was shown to be true for $k=r$ above.)

From \eqnref{fs} for $i>k$, we get when restricted to $W_k$
\begin{equation}
\rho(F_{i,m})=x_{ii,m}-\sum_{p\in\mathbb Z}\sum_{j=i+1}^nx_{ij,p}\partial_{x_{i+1,j,p-m}}
\end{equation}
and as a consequence
\begin{align*}
\rho(F_{k+1,m}^p)(u_k\otimes v)=x_{k+1,k+1,m}^pu_k\otimes v=u\otimes v
\end{align*}
for any $m\in\mathbb Z$, $p\in\mathbb N$. Since we assume $u_k\otimes v\in U(\widehat{\mathfrak g}) W_{\mathfrak a}(\tilde \lambda_r)$, the above equation tells us that $u\otimes v\in U(\widehat{\mathfrak g}) W_{\mathfrak a}(\tilde \lambda_r)$.

Now recall from the proof of \thmref{realization} (see \cite{MR2065632})
\begin{align}\label{serref}[\rho(F_i)&(z),\rho(F_j)(w)]
=\left(\delta_{i,j+1}a_{j,j+1}(w)-\delta_{j,i+1}a_{i,i+1}(z)
   \right)\delta(z-w)\notag \\
       &\quad
+\left(\delta_{i,j+1}\sum_{q=j+2}^{n}a_{jq}(z)a_{j+2,q}^*(z)-
   \delta_{j,i+1}\sum_{q=i+2}^{n}a_{iq}(z)a_{i+2,q}^*(z)
\right)\delta(z-w).
\end{align}
We can show by induction that for $1\leq j<i$,
\begin{align}\label{rootvector}
[\rho(F_i)(z_i),\dots, \rho(F_j)(z_j)]
&=\left(a_{ji}(z_j)+\sum_{q=i+1}^na_{jq}(z_j)a_{i+1,q}^*(z_j)
\right)\prod_{q=j}^{i-1}\delta(z_q-z_{q+1}).
\end{align}
Indeed this is true for $j=i-1$ by  \eqnref{serref} and the induction step is given by
\begin{align*}
[[\rho(F_i)(z_i)&,\dots, \rho(F_{j+1})(z_{j+1})], F_j(z_j)] \\
&=[a_{j+1,i}(z_{j+1})+\sum_{q=i+1}^na_{j+1,q}(z_{j+1})a_{i+1,q}^*(z_{j+1}), a_{jj}(z_j) \\
&\quad +\sum_{l=j+1}^na_{jl}(z_j)a^*_{j+1,l}(z_j)] \times \prod_{q=j+1}^{i-1}\delta(z_q-z_{q+1})\\
&=\left(a_{ji}(z_{j})+\sum_{q=i+1}^na_{jq}(z_{j})a_{i+1,q}^*(z_{j})\right)
\prod_{q=j}^{i-1}\delta(z_q-z_{q+1}).
\end{align*}
(There are no double contractions occurring in the calculation.)

Suppose we have shown $u^k_{i,\mathbf m}u_k\otimes v\in U(\widehat{\mathfrak g}) W_{\mathfrak a}(\tilde \lambda_r)$ for all $u^k_{i,\mathbf m}$ of the form $u^k_{i,\mathbf m}=x_{i,k+1,m_i}^{p_i}\cdots x_{k+1,k+1,m_{k+1}}^{p_{k+1}}$ and $u_k\in  \mathbb C_k[\mathbf x]$, $ m_j\in \mathbb Z$, $p_j\in\mathbb N
$  and $v\in \mathbb C_n[\mathbf y]$.  Then by \eqnref{rootvector} we get for $j<k+1$
\begin{align}
[\rho(F_{k+1})(z_{k+1}),\dots, \rho(F_j)(z_j)]_m
&=a_{j,k+1,m}+\sum_{q=k+2}^n\sum_{p\in \mathbb Z}a_{jq,p}a_{k+2,q,m-p}^*,
\end{align}
and thus
\begin{align*}
u^k_{i-1,\mathbf m}u_k\otimes v&= x_{i-1,k+1,m_{i-1}}^{p_{i-1}}x_{i,k+1,m_i}^{p_i}\cdots x_{k+1,k+1,m_{k+1}}^{p_{k+1}}u_k\otimes v \\
&=[\rho(F_{k+1})(z_{k+1}),\dots, \rho(F_{i-1})(z_{i-1})]_{m_{i-1}}^{p_{i-1}}(u_k\otimes v) \\
&=u^k_{i-1,\mathbf m}u_k\otimes v=u\otimes v
\end{align*}

This proves that $W_{k+1}\subset U(\widehat{\mathfrak
g})W_{\mathfrak a}(\tilde\lambda_r)$.  Hence by induction we have
$W_{\mathfrak g}(\lambda)=\mathbb C[\mathbf x]\otimes \mathbb
C_n[\mathbf y]=U(\widehat{\mathfrak g})W_{\mathfrak
a}(\tilde\lambda_r)$.
\end{proof}

We immediately have

\begin{cor}\label{irred}  If $\tilde\lambda\in\mathcal H^*$ is generic so
that $M_{r}(\tilde{\lambda})$ is irreducible, then $$M_{r}(\tilde
\lambda)\cong W_{\mathfrak g}(\tilde\lambda).$$
\end{cor}

\begin{proof}
For a generic $\tilde\lambda$, $W_{\mathfrak
a}(\tilde\lambda_r)\simeq W$ and $W_{\mathfrak
g}(\tilde\lambda)\simeq \tilde{W}$. It remains to apply
Corollary~\ref{cor20}.
\end{proof}

\section{Submodule structure of  Intermediate Wakimoto modules}

We saw in the previous section that the Intermediate Wakimoto
module $W_{\mathfrak g}(\tilde{\lambda})$ is generated by its
$\mathfrak a$-submodule $W_{\mathfrak a}(\tilde\lambda_r)$. We
will show in this section that in the generic case $W_{\mathfrak
a}(\tilde\lambda_r)$ determines completely the structure of
$W_{\mathfrak g}(\tilde{\lambda})$.

We assume in this case section that the Intermediate Wakimoto
module $W_{\mathfrak g}(\tilde{\lambda})$ is in general position,
namely that $\mathfrak B$ is non-degenerate. Hence, in particular,
$\tilde{\lambda}(c)\neq 0$. Note that it does not imply that
$W_{\mathfrak g}(\tilde{\lambda})$ is isomorphic to a
corresponding Verma type module, whose structure is known by
Theorem~\ref{vermatype}.

Consider a parabolic subalgebra $$\mathfrak p={B}_r+\hat{\mathfrak
a}$$ of $\hat{\mathfrak g}$ with the Levi factor $\tilde{\mathfrak
a}=\hat{\mathfrak a}+\hat{\mathfrak H}$ and the radical $\mathfrak
R$. Then $W_{\mathfrak a}(\tilde\lambda_r)$ belongs to the
standard category $\mathcal O$ of $\tilde{\mathfrak a}$-modules.
Moreover, we immediately have

\begin{lem}\label{le20}
$W_{\mathfrak a}(\tilde\lambda_r)$ is a $\mathfrak p$-submodule of
$W_{\mathfrak g}(\tilde{\lambda})$ with a trivial action of the
radical $\mathfrak R$.

\end{lem}

Hence we can consider $W_{\mathfrak a}(\tilde\lambda_r)$ as a
$\mathfrak p$-module with a trivial action of $\mathfrak R$ and
construct a generalized Verma module
$$M(W_{\mathfrak a}(\tilde\lambda_r)=U(\hat{\mathfrak g})\otimes_{U(\mathfrak p)}W_{\mathfrak a}(\tilde\lambda_r).$$

Lemma~\ref{le20} and Theorem~\ref{thm40} immediately imply that
there exists a canonical epimorphism $$\phi:M(W_{\mathfrak
a}(\tilde\lambda_r))\rightarrow W_{\mathfrak g}(\tilde{\lambda}).
$$
Hence, the module $W_{\mathfrak g}(\tilde{\lambda})$ is a
homomorphic image of the generalized Verma module $M(W_{\mathfrak
a}(\tilde\lambda_r))$.

Denote $M^f=1\otimes W_{\mathfrak a}(\tilde\lambda_r)$. The
following result describes the structure of the module
$M(W_{\mathfrak a}(\tilde\lambda_r))$.

\begin{thm}\label{thm50}
Let $N\neq 0$ be a submodule of $M(W_{\mathfrak
a}(\tilde\lambda_r))$ and $N^f=N\cap M^f$.
\begin{itemize}
\item[(i)] $N^f\neq 0$; \item[(ii)] If $\mathfrak B$ is
non-degenerate then $N\simeq U(\mathfrak g)\otimes_{U(\mathfrak
p)}N^f$, where $\mathfrak R$ acts trivially on $N^f$.
\end{itemize}
\end{thm}

\begin{proof}
The first statement follows from Lemma~17 in \cite{MR1858493}, while the
second statement follows from Theorem~8,(i) in \cite{MR1858493}.

\end{proof}

Using Theorem~\ref{thm50} we obtain the following description of
the submodule structure of the Intermediate Wakimoto modules in
the generic case.

\begin{cor}\label{cor30}
 Let $N\neq 0$ be a submodule of $W_{\mathfrak
g}(\tilde{\lambda})$.
\begin{itemize}
\item[(i)] $N\cap W_{\mathfrak a}(\tilde\lambda_r)\neq 0$;
\item[(ii)] If $\mathfrak B$ is non-degenerate then $N\simeq
U(\mathfrak g)\otimes_{U(\mathfrak p)}(N\cap W_{\mathfrak
a}(\tilde\lambda_r))$.
\end{itemize}

\end{cor}

\section{Conclusion}

As it was mentioned above Wakimoto modules can be obtained from
the classical Verma modules by an infinite number of twistings.
The same twisting can be applied to a Verma type module
$M_r(\tilde{\lambda})$ in the part $M(\tilde\lambda_r)$, i.e. only
using the reflections corresponding to the roots of $\mathfrak a$.
We will say that in this case the module is obtained by {\em real
twisting}. Clearly, imaginary Verma modules do not admit any real
twisting, while on the other hand any intermediate Wakimoto module
is obtained from the corresponding Verma type module by an an
infinite number of real twistings. Hence, all boson type
realizations associated with the natural Borel subalgebra
correspond to infinite (or empty in the imaginary case) real
twistings of corresponding Verma type modules. But we do not get
realizations of Verma type modules this way. In order to construct
boson type realizations for these modules one needs to start with
the Borel subalgebra different from the standard one or the
natural one and consider a corresponding "flag manifold". Its
 "cells" will produce boson type realizations for Verma type
modules, their contragradient analogs and finite real twistings.
We are going to address this question in a subsequent paper.

\section{Appendix}

In this section we present the formulas for the singular elements
in Imaginary Verma modules recently obtained by B.Wilson
\cite{Wi}. These formulas were inspired by the free field
realization of Imaginary Verma modules for $\hat{sl(2)}$.

Let $e_i$, $f_i$, $h_i$ be  a basis of $\hat{\mathfrak{sl}(2)}=\mathfrak{sl}(2)\otimes
\C[t, t^{-1}]$, $i\in \Z$ where $x_i:=x\otimes t^i$, $x\in\mathfrak{sl}(2)$. Consider the imaginary Verma module
$M_0(0)$ with a trivial highest weight. Then it has a submodule
generated by $h_i\otimes 1, i<0$. Denote by $V(0)$ the
corresponding quotient. A nonzero element $v\in V(0)$ is called
{\it singular} if $e_i v=0$ and $h_jv=0$ for all $i\in\mathbb Z$ and $j\in\mathbb N^*$. Let $F=\sum_i \C f_i$.
Then $V(0)$ is a free $U(F)$-module, where $U(F)$ is just a
polynomial algebra and the elements
$$
\Pi_{i=0}^{r-1}f_{s_i},\,\,\, s_i\in \Z,\,\,\, s_i\leq s_{i+1}, r\geq 1
$$
form a basis of $V(0)$.

Denote by $M_r=span_{\C}\{\Pi_{i=0}^{r-1}f_{s_i}| s_i\in \Z\}$.
Let $S_r$ be the set of singular vectors in $M_r$ and let $Sym(r)$
be the symmetric group realized as permutations of $0, 1, \ldots,
r-1$. For any $r\geq 1$ and $\bar{s}\in \Z^r$ set
$$v_r(\bar{s})=\sum_{\sigma\in Sym(r)}(-1)^{sgn \sigma}f_{s_0+\sigma(0)}f_{s_1+\sigma(1)}\ldots f_{s_{r-1}+\sigma(r-1)}.$$

\begin{thm}[\cite{Wi}]
Elements $v_r(\bar{s})$, $r\geq 1$, $\bar{s}\in \Z^r$, form a
basis of $S_r $.

\end{thm}

\bibliographystyle{alpha}
\bibliography{math}

\begin{thebibliography}{FKM01}

\bibitem[AL03]{MR1985191}
H.~H. Andersen and N.~Lauritzen.
\newblock Twisted {V}erma modules.
\newblock In {\em Studies in memory of Issai Schur (Chevaleret/Rehovot, 2000)},
  volume 210 of {\em Progr. Math.}, pages 1--26. Birkh\"auser Boston, Boston,
  MA, 2003.

\bibitem[BC94]{MR95i:17007}
Stephen Berman and Ben Cox.
\newblock Enveloping algebras and representations of toroidal {L}ie algebras.
\newblock {\em Pacific J. Math.}, 165(2):239--267, 1994.

\bibitem[BF90]{BF}
D.~Bernard and G.~Felder.
\newblock Fock representations and {B}{R}{S}{T} cohomology in ${\rm {s}{l}}(2)$
  current algebra.
\newblock {\em Comm. Math. Phys.}, 127(1):145--168, 1990.

\bibitem[BGG73]{MR55:2941}
I.~N. Bern{\v{s}}te{\u\i}n, I.~M. Gel{\cprime}fand, and S.~I. Gel{\cprime}fand.
\newblock Schubert cells, and the cohomology of the spaces ${G}/{P}$.
\newblock {\em Uspehi Mat. Nauk}, 28(3(171)):3--26, 1973.

\bibitem[BS83]{MR85g:81096}
N.~N. Bogoliubov and D.~V. Shirkov.
\newblock {\em Quantum fields}.
\newblock Benjamin/Cummings Publishing Co. Inc. Advanced Book Program, Reading,
  MA, 1983.
\newblock Translated from the Russian by D. B. Pontecorvo.

\bibitem[CF01]{MR2001k:17011}
B.~Cox and V.~Futorny.
\newblock Borel subalgebras and categories of highest weight modules for
  toroidal {L}ie algebras.
\newblock {\em J. Algebra}, 236(1):1--28, 2001.

\bibitem[CF04]{MR2065632}
Ben~L. Cox and Vyacheslav Futorny.
\newblock Intermediate {W}akimoto modules for affine {$\germ s\germ l(n+1,\Bbb
  C)$}.
\newblock {\em J. Phys. A}, 37(21):5589--5603, 2004.

\bibitem[CM01]{MR2002f:17041}
Charles~H. Conley and Christiane Martin.
\newblock A family of irreducible representations of the {W}itt {L}ie algebra
  with infinite-dimensional weight spaces.
\newblock {\em Compositio Math.}, 128(2):153--175, 2001.

\bibitem[Con74]{MR0374214}
Nicole Conze.
\newblock Alg\`ebres d'op\'erateurs diff\'erentiels et quotients des alg\`ebres
  enveloppantes.
\newblock {\em Bull. Soc. Math. France}, 102:379--415, 1974.

\bibitem[Cox94a]{MR95d:17026}
Ben Cox.
\newblock Structure of the nonstandard category of highest weight modules.
\newblock In {\em Modern trends in Lie algebra representation theory (Kingston,
  ON, 1993)}, volume~94 of {\em Queen's Papers in Pure and Appl. Math.}, pages
  35--47. Queen's Univ., Kingston, ON, 1994.

\bibitem[Cox94b]{C}
Ben Cox.
\newblock Verma modules induced from nonstandard {B}orel subalgebras.
\newblock {\em Pacific J. Math.}, 165(2):269--294, 1994.

\bibitem[Cox02]{MR2003g:17034}
Ben~L. Cox.
\newblock Two realizations of toroidal {$\germ s\germ l\sb 2(\Bbb C)$}.
\newblock In {\em Recent developments in infinite-dimensional Lie algebras and
  conformal field theory (Charlottesville, VA, 2000)}, volume 297 of {\em
  Contemp. Math.}, pages 47--68. Amer. Math. Soc., Providence, RI, 2002.

\bibitem[Cox04]{C2}
Ben Cox.
\newblock Fock space realizations of imaginary {V}erma modules.
\newblock {\em Algebras and Their Representations}, 8(2):173 -- 206, 2004.

\bibitem[CP87]{MR88h:17022}
Vyjayanthi Chari and Andrew Pressley.
\newblock A new family of irreducible, integrable modules for affine {L}ie
  algebras.
\newblock {\em Math. Ann.}, 277(3):543--562, 1987.

\bibitem[dBF97]{MR1482938}
Jan de~Boer and L{\'a}szl{\'o} Feh{\'e}r.
\newblock Wakimoto realizations of current algebras: an explicit construction.
\newblock {\em Comm. Math. Phys.}, 189(3):759--793, 1997.

\bibitem[DFO94]{MR1308982}
Yu.~A. Drozd, V.~M. Futorny, and S.~A. Ovsienko.
\newblock Harish-{C}handra subalgebras and {G}el\cprime fand-{Z}etlin modules.
\newblock In {\em Finite-dimensional algebras and related topics (Ottawa, ON,
  1992)}, volume 424 of {\em NATO Adv. Sci. Inst. Ser. C Math. Phys. Sci.},
  pages 79--93. Kluwer Acad. Publ., Dordrecht, 1994.

\bibitem[FBZ01]{MR1849359}
Edward Frenkel and David Ben-Zvi.
\newblock {\em Vertex algebras and algebraic curves}, volume~88 of {\em
  Mathematical Surveys and Monographs}.
\newblock American Mathematical Society, Providence, RI, 2001.

\bibitem[Fe{\u\i}84]{MR740035}
B.~L. Fe{\u\i}gin.
\newblock Semi-infinite homology of {L}ie, {K}ac-{M}oody and {V}irasoro
  algebras.
\newblock {\em Uspekhi Mat. Nauk}, 39(2(236)):195--196, 1984.

\bibitem[FF88]{MR89k:17016}
B.~L. Fe{\u\i}gin and {\`E}.~V. Frenkel.
\newblock A family of representations of affine {L}ie algebras.
\newblock {\em Uspekhi Mat. Nauk}, 43(5(263)):227--228, 1988.

\bibitem[FF90a]{MR92f:17026}
Boris~L. Fe{\u\i}gin and Edward~V. Frenkel.
\newblock Affine {K}ac-{M}oody algebras and semi-infinite flag manifolds.
\newblock {\em Comm. Math. Phys.}, 128(1):161--189, 1990.

\bibitem[FF90b]{MR92d:17025}
Boris~L. Feigin and Edward~V. Frenkel.
\newblock Representations of affine {K}ac-{M}oody algebras and bosonization.
\newblock In {\em Physics and mathematics of strings}, pages 271--316. World
  Sci. Publishing, Teaneck, NJ, 1990.

\bibitem[FKM01]{MR1858493}
Vyacheslav Futorny, Steffen K{\"o}nig, and Volodymyr Mazorchuk.
\newblock Categories of induced modules for {L}ie algebras with triangular
  decomposition.
\newblock {\em Forum Math.}, 13(5):641--661, 2001.

\bibitem[FS93]{FS}
Viatcheslav Futorny and Halip Saifi.
\newblock Modules of {V}erma type and new irreducible representations for
  affine {L}ie algebras.
\newblock In {\em Representations of algebras (Ottawa, ON, 1992)}, pages
  185--191. Amer. Math. Soc., Providence, RI, 1993.

\bibitem[Fut94]{MR95a:17030}
V.~M. Futorny.
\newblock Imaginary {V}erma modules for affine {L}ie algebras.
\newblock {\em Canad. Math. Bull.}, 37(2):213--218, 1994.

\bibitem[Fut97]{MR1627814}
Vyacheslav~M. Futorny.
\newblock {\em Representations of affine {L}ie algebras}, volume 106 of {\em
  Queen's Papers in Pure and Applied Mathematics}.
\newblock Queen's University, Kingston, ON, 1997.

\bibitem[Hua98]{MR99m:81001}
Kerson Huang.
\newblock {\em Quantum field theory}.
\newblock John Wiley \& Sons Inc., New York, 1998.
\newblock From operators to path integrals.

\bibitem[JK85]{JK}
H.~P. Jakobsen and V.~G. Kac.
\newblock A new class of unitarizable highest weight representations of
  infinite-dimensional {L}ie algebras.
\newblock In {\em Nonlinear equations in classical and quantum field theory
  (Meudon/Paris, 1983/1984)}, pages 1--20. Springer, Berlin, 1985.

\bibitem[JK89]{MR89m:17032}
Hans~Plesner Jakobsen and Victor Kac.
\newblock A new class of unitarizable highest weight representations of
  infinite-dimensional {L}ie algebras. {I}{I}.
\newblock {\em J. Funct. Anal.}, 82(1):69--90, 1989.

\bibitem[Kac90]{K}
Victor~G. Kac.
\newblock {\em Infinite-dimensional {L}ie algebras}.
\newblock Cambridge University Press, Cambridge, third edition, 1990.

\bibitem[Kac98]{MR99f:17033}
Victor Kac.
\newblock {\em Vertex algebras for beginners}.
\newblock American Mathematical Society, Providence, RI, second edition, 1998.

\bibitem[Kho]{Kh}
A.~Khomenko.
\newblock Some applications of gelfand-tsetlin modules.

\bibitem[MN99]{MR2000k:17036}
Atsushi Matsuo and Kiyokazu Nagatomo.
\newblock {\em Axioms for a vertex algebra and the locality of quantum fields},
  volume~4 of {\em MSJ Memoirs}.
\newblock Mathematical Society of Japan, Tokyo, 1999.

\bibitem[PRY96]{MR1412381}
J.~L. Petersen, J.~Rasmussen, and M.~Yu.
\newblock Free field realization of {${\rm SL}(2)$} correlators for admissible
  representations, and {H}amiltonian reduction for correlators.
\newblock {\em Nuclear Phys. B Proc. Suppl.}, 49:27--34, 1996.
\newblock Theory of elementary particles (Buckow, 1995).

\bibitem[Vor93]{MR1223226}
Alexander~A. Voronov.
\newblock Semi-infinite homological algebra.
\newblock {\em Invent. Math.}, 113(1):103--146, 1993.

\bibitem[Vor99]{MR1713301}
Alexander~A. Voronov.
\newblock Semi-infinite induction and {W}akimoto modules.
\newblock {\em Amer. J. Math.}, 121(5):1079--1094, 1999.

\bibitem[Wak86]{W}
Minoru Wakimoto.
\newblock Fock representations of the affine {L}ie algebra ${A}\sp {(1)}\sb 1$.
\newblock {\em Comm. Math. Phys.}, 104(4):605--609, 1986.

\bibitem[Wil05]{Wi}
B.~Wilson.
\newblock Structure of imaginary verma modules for affine lie algebras.
\newblock 2005.

\end{thebibliography}
\def\cprime{$'$}

\end{document}